\documentclass[12pt]{article}

\usepackage{apacite,latexsym,amsmath,amssymb,amsthm}
\usepackage{natbib}

\begin{document}

\title{Directionally collapsible parameterizations of multivariate binary distributions}
\author{Tam\'as Rudas\\Department of Statistics, Faculty of Social
Sciences\\
E\"{o}tv\"{o}s Lor\'{a}nd University, Budapest\\
\texttt{rudas@tarki.hu}}
\maketitle

\begin{abstract} \noindent Odds ratios and log-linear parameters are not collapsible, meaning that  including a variable into the analysis or omitting one from it, may change the strength of association among the remaining variables. Even the direction of association may be reversed, a fact that is often discussed under the name of Simpson's paradox. A parameter of association is directionally collapsible, if this reversal cannot occur. The paper investigates the existence of parameters of association which are directionally collapsible. It is shown, that subject to two simple assumptions, no parameter of association, which depends only on the conditional distributions, like the odds ratio does, can be directionally collapsible. The main result is that every directionally collapsible parameter of association gives the same direction of association as a linear contrast of the cell probabilities does. The implication for dealing with Simpson's paradox is that there is exactly one way to associate direction with the association in any table, so that the paradox never occurs.\end{abstract}

\hspace{5mm}\small{Keywords: directional collapsibility, odds ratio, Simpson's paradox, parameterization

\hspace{4mm} of binary distributions, variation independence from marginal distributions}

\section{Introduction}

This paper studies the relationships between certain properties that parameters of associations for binary distribution may have. \cite{GK} gave an overview of bivariate parameters of association and they argued that no single concept of association may be used in all research problems. Interest since then has turned towards the multivariate case and, although there have been alternative suggestions, see, e.g., \cite{bah61}, \cite{lancaster69}, applications and theoretical work in the last fifty years have concentrated around the odds ratio and quantities derived from it, mostly because of their relevance in log-linear and other graphical Markov models, see, e.g, \cite{bfh75}, \cite{lauritzen96}. The multivariate version of the odds ratio was first considered in \cite{ba35}, see also \cite{bi63},  and \cite{kk68}  for a review of related approaches. However, not every analyst is entirely satisfied with odds ratios (or their logarithms), as parameters of association. First, the standard error of the sample odds ratio, as an estimator, depends not only on the true value of the odds ratio, but is a monotone function of the sum of the reciprocals of the cell probabilities, resulting in high variability of estimators. Second, lack of collapsibility is often cited as an undesirable property, see, e.g., \cite{wh78}, \cite{we87} and \cite{ve14}. The fact that even the direction of association may change after collapsing (e.g.,  taking the new drug may be associated with recovery for both male and female patients, but disregarding sex, taking the old drug is associated with recovery) is seen as paradoxical by many, as shown by the widespread literature on 'Simpson's paradox'. In addition to well-known occurrences of Simpson's paradox in sociology, education and the health sciences, it is being discussed  in genetics (\citealp{br14}) and in  physics (\citealp{li13}).

 As opposed to the vast majority of this literature, Simpson's paradox is not considered here as a special, perhaps negative, feature of the data for which it occurs, rather it is considered as a characteristic of the parameter of association applied, namely the odds ratio, that conditional and marginal associations may have opposing directions (cf. \citealp{we87, rudas10}). Directional collapsibility means, that such a reversal cannot occur.

The direction of association is readily interpreted for $k=2$. If one variable is treatment, the other is response to it, then the direction of association tells whether the treatment is beneficial or detrimental to the response.  If the two variables are treated on an equal footing, that is, none of them is assumed to be a response to the other, then the direction of association tells whether concordant or discordant types of observations are more likely. For more than $2$ variables, when one is response to the others, if all treatments are beneficial when applied individually, the direction of association may tell   whether applying all treatments has additional benefit, or whether it is beneficial, at all. However, just like there is no single parameter of association, there is also no single meaning of association. When the variables are treated on an equal footing, one possible interpretation is given in (\ref{recursion}) and in the discussion following it.

The paper investigates the possibility of finding directionally collapsible parameters of association, which also provide a parameterizations of multivariate binary distributions. The main results are obtained under two simple assumptions made for parameters of association, which are described and motivated in Section 2. These two properties are possessed not only by the odds ratio, but also by a simple contrast of the cell probabilities defined in (\ref{DI}). It is also shown in Section 2, that both the odds ratios and the contrasts, associated with all marginal distributions, constitute a parameterization of the joint distribution.

The main results of the paper are given in Section 3. Variation independence of the odds ratio from lower dimensional marginal distributions, formulated here as dependence on the conditional distributions only, which is a very desirable property in other contexts (see, e..g, \citealp{rudas98}),  turns out to imply the lack of directional collapsibility. More precisely, any parameter of association which depends only on the conditional distributions,  assigns the same direction of association to every distribution as the log odds ratio does, and, therefore, is not directionally collapsible. On the other hand, a parameter of association is directionally collapsible, if and only if it assigns the same direction of association to every distribution, as the contrast of the cell probabilities does.   

One is then left with the following simple situation. If the two properties described in Section 2 are assumed, then all parameters of association which depend on the conditional distributions only, deem the direction of association as the odds ratio does, and are not directionally collapsible. Further, all directionally collapsible parameters of association assign the same direction to association as the contrast does, and the latter also provides a parameterization of the distribution.

Section 4 concludes the paper with a brief discussion of the potential use of the contrast as a parameter of association and the implications for dealing with Simpson's paradox. Those analysts who are interested in the direction of association only, and find the contrast being overly simple, failing to properly describe their concept of association, cannot avoid Simpson's paradox and  have to learn to accept the reversal as not paradoxical. On the other hand, those, whose main concern remains to avoid Simpson's paradox, and are ready to use the contrast to determine the direction of association, will be happy to see that the contrast has very attractive sampling properties, including that its sampling distribution does not depend on the number of variables involved, rather only on its population value.

\section{Some properties of parameters of association}

This paper deals with parameters of the joint distribution of $k$ binary variables. Such distributions may be written as entries in the cells $t$ of a $2^k$ contingency table, $T_k$. The cells of such a table may be identified with sequences of $1$'s and $2$'s of length $k$, and the notation $t = (j_1, j_2, \ldots, j_k)$ will be used, where $j_i$ is $1$ or $2$ for all $ i = 1, 2, \ldots, k $. The distributions to be considered are not restricted to probability distributions summing to $1$ and not even to frequency distributions with integer values. The set of any positive entries $(p(t), \, t \in T_k)$ in the contingency table$T_k$  will be called a distribution. 

This paper offers no definition of what is a parameter of association, rather the relationships between different possible characteristics are investigated. As pointed out by \citet{GK}, see also \cite{da74}, \cite{stre90}, \cite{rb04}, there are several ways to define parameters of association, which may be relevant in different research context, and these different parameters may have different characteristics. The properties which will be assumed here, seem appropriate in the common situations when:

(i) The variables considered describe the presence (category 1) or absence (category 2) of various characteristics. Association means that these characteristics tend to occur together. If the joint distribution is uniform, there is no association, and the stronger is the tendency for all the characteristics to occur together, the stronger is the association.

(ii) Association has a direction, and any pattern of association of $k-1$ variables, combined with the presence or the absence of the $k$-th characteristic imply different directions of association.

These assumptions are formulated as Properties 1 and 2.

\textbf{Property 1.}
A parameter of association $f_k$  is  a continuous real function on the set of distributions $(p(t): \,\, t \in T_k)  $, such that 
\begin{equation}\label{zero}
f_k(p(t): \,\, t \in T_k)  = 0 \makebox{ if } p(t)=c \makebox{ for all } t \in T_k  
\end{equation}
and $f_k$ is strictly monotone  increasing in $p(1, 1, \ldots, 1)$.

\hfill{} \qed

\textbf{Property 2.}
If $(p(t), \, t \in T_k)$ and   $(q(t), \, t \in T_k)$  are  distributions, such that there is an $i \in \{1, 2, \ldots, k\}$ with 
$$
p(j_1, \ldots, j_{i-1}, j_i, j_{i+1}, \ldots, j_k) = q(j_1, \ldots, j_{i-1}, j_i^*, j_{i+1}, \ldots, j_k), 
$$  
where $j_i^*+j_i=3$, for all $(j_1, \ldots, j_{i-1},  j_{i+1}, \ldots, j_k)$, then
$$
sgn \left(  f_k(p(t): t \in T_k) \right) = - sgn \left(  f_k(q(t): t \in T_k) \right).
$$
\hfill{} \qed

This is not the most parsimonious formulation of these assumptions: Property 2 implies (\ref{zero}). Let the cells  of $T_k$ with an even number of $2$'s be denoted by   $T_{ke}$ and those with an odd number of $2$'s by $T_{ko}$. Swapping the categories of a variable, as described in Property 2, interchanges these two subsets of the cells.

Interaction parameters which are contrasts between certain functions of the cell entries play a central role in this paper. More precisely, let $h$ be a  monotone increasing continuous real function and  consider

\begin{equation}\label{contrast}
f_k(p(t): t \in T_k) =  \sum_{t \in T_{ke}} h(p(t)) - \sum_{t \in T_{ko}} h(p(t)) .
\end{equation}
Because $(1,1, \ldots, 1) \in T_{ke}$, Property 1 holds, and because if
$$
(j_1, \ldots, j_{i-1}, j_i, j_{i+1}, \ldots, j_k) \in T_{ke},
$$
then 
$$
(j_1, \ldots, j_{i-1}, j_i^*, j_{i+1}, \ldots, j_k) \in T_{ko},
$$
and vice versa, Property 2 holds, too, for interaction parameters of the type (\ref{contrast}).

If $f_k$ is of the form (\ref{contrast}), then it may be written as
\begin{equation}\label{di}
f_k(p(t): t \in T_k) = \sum_{t \in T_k} (-1)^{e't-k} h(p(t)),
\end{equation}
where $e'$ is the transpose of a column vector of length $k$, consisting of $1'$s. 

The following example illustrates parameters of association of the type (\ref{contrast}). 

\textbf{Example 1.} 
The $k-1$st order odds ratio for a $k$-dimensional distribution is 
\begin{equation}\label{OR}
OR_k(p(t): t \in T_k) = \frac{\prod_{t \in T_{ke}} p(t)}{\prod_{t \in T_{ko}} p(t)}
\end{equation}
and log $OR_k$ is an interaction parameter, to be denoted as $LOR_k$. The log odds ratios are closely related to the log-linear parameters of the distribution (see, e.g., \citealp{rudas98}). 

The log odds ratios may also be generated as in (\ref{contrast}), by using $h=\log$:
$$
LOR_k(p(t): t \in T_k) =  \sum_{t \in T_{ke}} \log(p(t)) - \sum_{t \in T_{ko}} \log(p(t)) .
$$

The difference parameter of association is
\begin{equation}\label{DI}
DI_k(p(t): t \in T_k) = \sum_{t \in T_{ke}} p(t) - \sum_{t \in T_{ko}} p(t) ,
\end{equation}
which is obtained from (\ref{contrast}) by choosing $h$ as the identity function.

Finally, by choosing $h = \exp$ in (\ref{contrast}) gives 
$$
EX_k(p(t): t \in T_k) =  \sum_{t \in T_{ke}} \exp(p(t)) - \sum_{t \in T_{ko}} \exp(p(t)) .
$$
\hfill{} \qed

Parameters of association of the type (\ref{contrast}) are not only contrasts between functions of the entries in $ T_{ke} $ and  $ T_{ko} $, but also a comparison of the strengths of association in parts of the table defined by specific indices of a variable. Let  $ T_{k-1}(V_i=1) $ be the part of the table where the $ i $th variable is $ 1 $,  and  $ T_{k-1}(V_i=2) $ be the part of the table where the $j$ith variable is $2$.  These are $k-1$-dimensional tables formed by the variables other than $V_i$. Then, if $f_k$ is of the type (\ref{contrast}),  it may be obtained by the following recursion, irrespective of the choice of $i$ :
$$
f_1(p(t): t \in T_1) = h(p(1)) - h(p(2))
$$
\begin{equation}\label{recursion}
f_k(p(t): t \in T_k) = 
\end{equation}
$$
f_{k-1}(p(t) : t \in T_{k-1}(V_i=1)) - f_{k-1}(p(t): t \in T_{k-1}(V_i=2)).
$$
To see that (\ref{contrast}) and (\ref{recursion}) give the same, only the signs of the quantities $h(p(t))$ need to be checked. 
For every $t \in T_k$, the the sign of $h(p(t))$ in $f_k$ in (\ref{contrast}) is the same as the sign in $f_{k-1}$ in (\ref{recursion}), if and only if $V_i=1$, and is the opposite when $V_i=2$, because the sign depends on the parity of the number of $2$'s among the indices. This reversal is introduced in (\ref{recursion}) by the negative sign of the second term. 

Formula  (\ref{contrast}) may seem counter-intuitive, even "wrong", as it suggests, as implied by Property 2, that large entries in cells with an odd number of $2$'s among their indices imply weak association. Formula (\ref{recursion}) shows, that (\ref{contrast}) is a comparison, showing, for any variable $V_i$, the amount by which association is stronger, when, in addition to all other characteristics, also the one indicated by $V_i$ is present ($j_i=1$), as opposed to when it is not ($j_i=2$).

However, there are functions of the cell entries which possess Properties 1 and 2, but cannot be written in the form of (\ref{contrast}), as illustrated next.

\textbf{Example 2.} 
Let  let $d$ be strictly monotone but non-linear function.Then
$$
d (\sum_{t \in T_{ke}} p(t)) - d(\sum_{t \in T_{ko}} p(t)) .
$$
is a parameter of association which cannot be written in the form of (\ref{contrast}).

For example, in the case of $k=2$, with the usual notation,
$$
(p(1,1)+p(2,2))^3 - (p(1,2)+p(2,1))^3
$$
is not a linear contrast of any function of the cell entries.

\hfill{ } \qed

The next example illustrates parameters of association which do not possess Properties 1 and 2.
 
\textbf{Example 3.} 
One may say, that in the following distribution, the three variables (possessing the three characteristics) do show some association, because it is more likely to have all three characteristics present, than any other pattern of presence or absence.

$$
\begin{array}{|r|r|} \hline 
   0.3140 & 0.098    \\
\hline
    0.098 & 0.098 \\
\hline
\end{array}
\,\,\,\,
\begin{array}{|r|r|} \hline 
   0.098 & 0.098   \\
\hline
    0.098 & 0.098 \\
\hline
\end{array}
$$ 

Indeed, the Bahadur parameter \citep{bah61} associates the value of $0.103$ with this distribution. By the same argument, one might think that the association is stronger in the following distribution.

$$
\begin{array}{|r|r|} \hline 
   0.9965 & 0.0005    \\
\hline
    0.0005 & 0.0005 \\
\hline
\end{array}
\,\,\,\,
\begin{array}{|r|r|} \hline 
   0.0005 & 0.0005   \\
\hline
    0.0005 & 0.0005 \\
\hline
\end{array}
$$

However, the Bahadur parameter associates the value of  $-5.54$ with this distribution, indicating a negative association among the three variables, thus Property 1 does not hold. On the other hand, Property 1 does hold for the Bahadur parameter in the cae of $k=2$.

Parameters of association obtained by some normalization of the chi-squared statistic (see \citealp{GK}) are always nonnegative, thus cannot possess Property 2.

\hfill{ } \qed

\cite{rudas10} discussed treatment selection in the case of a single treatment and a single response variable. The conditions under which he showed that every decision rule which avoids Simpson's paradox for all data sets, chooses the same treatment as $DI_2$ does, are implied by Properties $1$ and $2$. 

An important property  of the interaction parameters $LOR_k$ and $DI_k$  is that they constitute a parameterization of the distributions on the contingency table.   Parameterization means that the vector valued function, which for every distribution on $T_k$ gives its $2^k$ interaction parameters (one for every  subset of the variables), is invertible.

For easier formulation of this fact, these interaction parameters are extended to apply to zero-dimensional subsets, so that $LOR_{0}$ is the logarithm of the product of the entries in the table, and $DI_{0}$ is their sum.

\textbf{Theorem 1.} Let $T_k$ be a  $k$-dimensional binary contingency table formed by the ranges of the variables $V_1, \ldots , V_k$. Let $m$ be a $0-1$ vector of length $k$,  and let $M$ be the set of all such vectors. Let $\mathcal{V}_m$ be the subset of the variables consisting of those $V_i$, for which $m_i$ is not zero. Finally, let all the  parameters of association 
\begin{equation}\label{212}
f_{{e'm}}(p(t): t \in T_{e'm}(\mathcal{V}_m)),\,\,\, m \in M,
\end{equation} 
where $e'm$ is the sum of the components of $m$ and $T_{e'm}(\mathcal{V}_m)$ is the contingency table with the joint distribution of the variables in $\mathcal{V}_m$, be given. Then, if $f_k=LOR_k$ or $f_k=DI_k$, the distribution  on $T_k$ may be reconstructed. 

\textbf{Proof.} In the case, when $f_k=LOR_k$, (\ref{212}) is essentially a marginal log-linear parameterizaton as described by \cite{br02}, with all subsets of the set of variables being a hierarchical and complete class, and the claim follows from their  Theorem 2, where a reconstruction algorithm based on repeated applications of the Iterative Proportional Fitting procedure was also described. 

In the case, when $f_k=DI_k$, the given interaction parameter values define a system of linear equations for the cell entries. To formulate the equations in this system, consider a vector $m$. Each entry in the marginal table defined by $m$, is the sum of those entries of $T_k$, which are in cells with such vectors of indices $t$, that are identical to each other in all the positions that have a $1$ in $m$. When $DI_{e'm}$ is computed for the marginal table defined by $m$, all these entries have the same sign, namely, the sign associated with the marginal entry in the $e'm$-dimensional table by $DI_{e'm}$, which is  
$$
(-1)^{t'm-1'm},
$$
as implied by (\ref{di}). Thus, the left hand side of  the equation associated with $m$ is
$$
\sum_{t \in T_k}     (-1)^{t' m- 1'm} p(t)    ,
$$
and the right hand side is the value of the parameter of association for the marginal defined by $m$. This system of equations does have a positive solution by assumption, and as  the $2^k \times 2^k$ matrix of coefficients is shown below to be of full rank, it only has one solution. 

To see the rank of the coefficient matrix, consider any two of its rows, say, the ones associated with  different vectors $m_1$ and $m_2$. There is a position, where one of these vectors is $1$, and the other one is $0$. To simplify notation, it is assumed that $m_1$ is $0$ and $m_2$ is $1$ in position $k$. Then, for any two cells that are identical in the first $k-1$ indices, but one has a $1$, the other has a $2$ in the $k$th position,  exactly $1$ will have identical signs in the two rows, and exactly $1$ will have different signs, because changing the last index from $1$ to $2$ leaves the sign of the entry in the row (i.e., equation) associated with $m_1$ unchanged, but changes the sign of the entry in the row (i.e., equation) associated with $m_2$, as the sign depends on the parity of the number of $2'$s  among the indices of the cells. Therefore, half of the entries have identical, and half of the entries have different signs in the two rows, thus the two rows of coefficients  are orthogonal. If any two rows of the coefficient matrix are orthogonal, then the matrix is of full rank.

Any algorithm to find the solution of a system of linear equations may be used to reconstruct the distribution in $T_k$.   
\hfill{ }\qed

\section{Directional collapsibility}

The central question in  this paper is directional collapsibility of parameters of association, which is now defined formally as Property 3.

\textbf{Property 3.} If for some $i \in \{1, \ldots, k \}$,
$$
sgn \left( f_{k-1}(p(t): t \in T_{k-1}(V_i=1))  \right) = sgn \left( f_{k-1}(p(t): t \in T_{k-1}(V_i=2))  \right),
$$
then also
$$
sgn \left( f_{k-1}(p(t): t \in T_{k-1}(V_i=+))  \right)
$$
$$
=sgn \left( f_{k-1}(p(t): t \in T_{k-1}(V_i=1)) \right) = sgn \left( f_{k-1}(p(t): t \in T_{k-1}(V_i=2))  \right),
$$
where $T_{k-1}(V_i=+)$ is obtained from $T_k$ by collapsing (marginalizing) over $V_i$ 

\hfil{ } \qed

\textbf{Example 4.} It is well known that $LOR_k$ is not directionally collapsible.  On the other hand, $DI_k$ is directionally collapsible. For simplicity of notation, this will be shown now for $i=1$. It follows form (\ref{di}), that, with $e$ being a vector of $1'$s of length $k-1$,
$$
DI_{k-1}(p(t): t \in T_{k-1}(V_1=j)) = \sum_{t_{k-1} \in T_{k-1}} (-1)^{e't_{k-1} - (k-1) }  p(j,t_{k-1}),
$$
where $t_{k-1}$ is a cell in $T_{k-1}$ and $(j,t_{k-1})$ is a cell in $T_k$. 
Then, with $(+,t_{k-1})$ being a marginal cell,
$$
DI_{k-1}(p(t): t \in T_{k-1}(V_1=+) ) = \sum_{t_{k-1} \in T_{k-1}} (-1)^{e't_{k-1} - (k-1)  }  p(+,t_{k-1})
$$
$$
= \sum_{t_{k-1} \in T_{k-1}} (-1)^{e't_{k-1} - (k-1) }  p(1,t_{k-1}) + (-1)^{e't_{k-1} - (k-1)  }  p(2,t_{k-1}) 
$$
$$
= DI_{k-1}(p(t): t \in T_{k-1}(V_1=1)) + DI_{k-1}(p(t): t \in T_{k-1}(V_1=2)),
$$
and then the sign of the left hand side is equal to the common sign of the terms on the right hand side, which is what was to be seen. In fact, the argument shows that $DI_k$ is not only directionally collapsible, but is also collapsible.

\hfil{ } \qed

The first result of this section identifies a property of the $LOR_k$, which implies its lack of directional collapsibility, and, consequently, all  parameters of association with this property also lack directional collapsibility. This property is that  the  value of the  parameter of association depends on the conditional distributions only,  in the  sense given in the next definition. 

\textbf{Property 4.} If the distributions $( p(t), t \in T_k) $ and  $ (q(t), t \in T_k) $ are such, that there exists a variable $V_i$, such that its conditional distributions, given the categories of all other variables, derived from $p$ and $q$,  coincide, then

$$
f_k(p(t): t \in T_k)= f_k(q(t): t \in T_k),
$$ 

\hfill{ } \qed

The condition for the equality of the conditional distributions, written for the first variable, is that
\begin{equation}\label{conddistr}
\frac{p(1, t_{k-1})}{p(+,t_{k-1})}  = \frac{q(1, t_{k-1})}{q(+,t_{k-1})},
\end{equation}
for all cells $t_{k-1}$ of the table formed by the ranges of the last $k-1$ variables.

A celebrated characteristic of the odds ratio is variation independence of $LOR_k$ from the marginal distribution of any $k-1$ variables. This  property is usually formulated (\citealp{rudas98}) by saying that if  $( r(t), t \in T_k) $ and  $ (s(t), t \in T_k) $ are distributions on $T_k$, then there always exists a distribution $( u(t), t \in T_k) $, that has the $k-1$ dimensional marginal distributions of the first distribution, and the $k-1$st order odds ratio of the second one. This form of definition is applied to avoid the problems stemming from the $k-1$ dimensional marginal distributions not being variation independent for $k>2$ among themselves, see \cite{br02}.  The theory of mixed parameterization of exponential families \citep{bn78} implies that there is only one  distribution $u$.  Property 4 implies this variation independence.

\textbf{Example 5.} Obviously, $LOR_k$ depends on the conditional distributions only but $DI_k$ does not have this property. 

\hfill{ }\qed

The next theorem shows that if Property 4 is assumed, then $f_k(p(t): \,\, t \in T_k)$ is  equal to the value of $f_k$ for a special distribution, derived from $p$.

\textbf{Theorem 2.} Let $f_k$ be a parameter of association with Property 4. Then, 
$$
f_k(p(t): \,\, t \in T_k) = f_k(q(t): \,\, t \in T_k),
$$
where the distribution  $(q(t): \,\, t \in T_k)$ is such, that 
$$
q(t) = 1 \makebox { if } t \neq (1, \ldots, 1) 
$$
$$
q(1, \ldots, 1) ) = OR_k(p(t): \,\,t \in T_k).
$$
The proof is based on a series of transformations, which are first illustrated for $k=3$.

\textbf{Example 6.} For $k=3$, write the distribution as follows:
$$
\begin{array}{|r|r|} \hline 
   p(111) & p(121)    \\
\hline
    p(211) & p(221) \\
\hline
\end{array}
\,\,\,\,
\begin{array}{|r|r|} \hline 
   p(112) & p(122)   \\
\hline
    p(212) & p(222) \\
\hline
\end{array}
$$
The first transformation is to divide both $p(1,j,k)$ and $p(2,j,k)$ by the latter, for all choices of $j$ and $k$, which yields
$$
\begin{array}{|r|r|} \hline 
   \frac{p(111)}{p(211)} & \frac{p(121)}{p(221)}    \\
\hline
    1 & 1 \\
\hline
\end{array}
\,\,\,\,
\begin{array}{|r|r|} \hline 
   \frac{p(112)}{p(212)} & \frac{p(122)}{p(222)}   \\
\hline
   1 & 1 \\
\hline
\end{array}
$$ 
The conditional distribution of $V_1$, given $V_2$ and $V_3$ in this distribution is the same as in $(p(t), \,\, t \in T_k)$, thus, if $f_3$ depends on the conditional distributions only, its value remains the same.

The next transformation is to divide the entry in cell $(1, 1, k)$ and in cell  $(1,2,k)$, for all choices of $k$, by the latter, yielding
$$
\begin{array}{|r|r|} \hline 
   \frac{p(111)}{p(211)} / \frac{p(121)}{p(221)} & 1    \\
\hline
    1 & 1 \\
\hline
\end{array}
\,\,\,\,
\begin{array}{|r|r|} \hline 
   \frac{p(112)}{p(212)} / \frac{p(122)}{p(222)} & 1  \\
\hline
   1 & 1 \\
\hline
\end{array}
$$ 
As the conditional distribution of the second variable, given the first and the third did not change, the value of $f_3$ is also unchanged.

The last  transformation is to divide the entries in cells $(1,1,1)$ and $(1,1,2)$ by the latter. 
This gives
in cell $(1,1,1)$
$$
\left( \frac{p(111)}{p(211)} / \frac{p(121)}{p(221)}\right)/\left( \frac{p(112)}{p(212)} / \frac{p(122)}{p(222)} \right) ,
$$
which is the $2$nd order odds ratio, and the other cells all contain $1$. The value of $f_3$ is still unchanged, as the last transformation left the conditional distribution of the third variable, given the first two, unchanged.
 
\hfill{ }\qed

\textbf{Proof.} The proof applies a series of transformations to the distribution in $T_k$, such that each step leads to a distribution with the same value of $f_k$ and at the end of the transformations, the odds ratio of the $k$ variables appears in cell $(1, \ldots, 1)$, and the other cells all contain $1$.

To define such a series of transformations, note that if both $p(1, t_{k-1})$ and $p(2, t_{k-1})$, for a fixed $t_{k-1}$, are multiplied by the same number, the value of $f_k$ remains unchanged, because the conditional distribution of the first variable, conditioned on the last $k-1$ variables, remains unchanged, and similarly for any variable other than $V_1$. There is one transformation for each variable and they are applied consecutively to the result of the previous transformation. 

The transformation for  variable $V_1$ is dividing both  $p(1, t_{k-1})$ and $p(2, t_{k-1})$,  by $p(2, t_{k-1})$, for all possible choices of the last $k-1$ indices.  The conditional distribution of $V_1$, given all other variables, does not change, so $f_k(p(t): t  \in T_k)$ will remain unchanged, too. The transformation will change the $p(2, t_{k-1})$ entries to $1$, and the new value of  $p(1, t_{k-1})$ will be
$$
p(1, t_{k-1}) / p(2, t_{k-1}).
$$

The next step is for variable $V_2$. It consists of dividing both $p(1, 1, t_{k-2})$ and $p(1, 2, t_{k-2})$ by the latter entries, for all choices of the $k-2$ indices in  $t_{k-2}$. Note that $p(2,t_{k-1})=1$, thus $p(2, 1, t_{k-2})$ and $p(2, 2, t_{k-2})$ are also equal to $1$ and need not be divided. This will leave the conditional distribution of $V_2$, given all other variables, unchanged, so the value of $f_k(p(t) : t \in T_k)$ is also unchanged. This transformation does not affect any cell of the form $(2, t_{k-1})$, thus the entry in any cell with $2$ as the first index remains $1$, and, in addition, the entries in every cell with a second index equal to $2$  also become equal to $1$. 

The $i$th step is applied to a table with entries in the cells having $2$ as any of their $1$st, $2nd$, ... , $i-1$st indices equal to $1$. It consist of dividing all cell entries  $p(1, \ldots, 1,1, t_{k-i})$ and $p(1, \ldots, 1,2, t_{k-i})$ by the latter. As the conditional distribution of $V_i$, given all other variables, remains unchanged, so does $f_k(p: p \in T_k)$, too.

The last step of the series of transformations is for variable $V_k$ and it consists of dividing the entries in cells $(1, \ldots, 1, 1)$ and $(1, \ldots, 1, 2)$ by the latter. This does not affect the value of the parameter of association and makes all cell entries, except for the one in $(1, \ldots, 1)$ equal to $1$.

All steps leave the value of $f_k(p(t): t \in T_k)$ unchanged. If an entry was made equal to $1$, it is not changed later during the transformations. After this series of transformations, all original cell entries appear in a multiplicative formula in cell $(1, 1,\ldots, 1)$, and all other entries are made equal to $1$, because these latter cells contain at least one index equal to $2$, and were divided by their own value in step $i$ of the transformation, if their first index equal to $2$ is the $i$th one.

Some of the original entries in $ (p(t): \,\, t \in T_k) $ appear in the numerator in cell  $(1,1, \ldots, 1)$, some appear in the denominator. The (original) value $p(1, 1, \ldots, 1)$ is in the numerator, because there is no division performed with that entry. All other terms appear in this cell as a result of a number of consecutive divisions. If during the series of transformations, a value $p(t_k)$ goes into the numerator of the entry in a cell, the next time, if such exists, when the entry in that cell is used for division, this value will appear in the denominator. Therefore, whether (the original)  $p(t_k)$ ends up in the numerator or in the denominator of the (final) entry in cell $(1, 1, \ldots, 1)$, depends on the parity of times, divisions involving that term occurred. And this is exactly the number of indices equal to $2$ in   $t_k$. If the number of $2$'s is even, the original value in the  cell  will be in the numerator, if it is odd, the original value will be in the denominator, which gives (\ref{OR}).

\hfill{ }\qed

Consequently, if Property 1 is also assumed, then the direction of association can be determined for parameters of association which depend on the conditional distributions only, as formulated in the next theorem.

\textbf{Theorem 3.} Let $f_k$ be a parameter of association with Properties 1 and 4. Then, 
$$
sgn \left( f_k(p(t): t \in T_k)  \right) = sgn \left(LOR_k(p(t) : t \in T_k)  \right) , 
$$
that is, $f_k$ assigns the same direction of association to all distributions as $LOR_k$ does, and, therefore,  $f_k$ is not directionally collapsible.

\textbf{Proof.} 
Consider  the distribution constructed in Theorem 2. If it had $1$ in every cell, then $f_k$ would be zero but it has the odds ratio of $(p(t): t \in T_k) $ in cell $(1, \ldots, 1)$.  Thus, to obtain this distribution from the one containing  $1$'s in every cell,  the entry in cell $(1, \ldots, 1)$ has to be increased / left unchanged / decreased, depending on whether the odds ratio is more than / equal to / less than $1$, making $f_k$ positive / zero / negative, which is also the sign of $LOR_k$.

The second claim of the theorem is implied by the first one. 

\hfill{ }\qed

Theorem 3 says that no parameter of association, which  depends on the conditional distributions only in the sense of Property 4,  can avoid Simpson's paradox for all data sets, if Property 1 is also assumed to hold. Every such parameter of association assigns the same sign to association to any distribution, as the $LOR_k$ does. Consequently, as long as one is only interested in the direction of association and wishes to use parameters of association  which  depend on the conditional distributions only, it is sufficient to use the $LOR_k$, but Simpson's paradox cannot be avoided. 

The next example illustrates that there are parameters of association, which do not depend on the conditional distributions only, yet are not directionally collapsible, thus the converse of Theorem 3 does not hold.

\textbf{Example 7.}
The parameter of association $EX_k$ does not depend on the conditional distributions only. In the following two distributions the conditional distribution of $V_1$ given $V_2$ is the same, yet the value of $EX_2$ for the first one is 79.67, and for the second one is -0.60.

$$
\begin{array}{|r|r|} 
\hline 
  2 & 3    \\
\hline
    4 & 5 \\
\hline
\end{array}
$$

$$
\begin{array}{|r|r|} \hline 
  0.6 & 0.6   \\
\hline
    1,2 & 1\\
\hline
\end{array}
$$

In spite of this, $EX_k$ is not directionally collapsible. In both of the following tables, $EX_2$ is positive (259.94 and 143.46, respectively)

$$
\begin{array}{|r|r|} 
\hline 
  6 & 5   \\
\hline
    3 & 3 \\
\hline
\end{array}
$$

$$
\begin{array}{|r|r|} \hline 
  5 & 7   \\
\hline
    1 & 7 \\
\hline
\end{array}
$$
but in the collapsed table

$$
\begin{array}{|r|r|} 
\hline 
  11 & 12   \\
\hline
    4 & 10 \\
\hline
\end{array}
$$
it is negative (-77694.70). Note that $LOR_k$ does not exhibit Simpson's paradox for these data.

\hfill{ }\qed

The main result of the section is that all directionally collapsible parameters of association judge the direction of association like $DI_k$ does, if Properties 1 and 2 are assumed to hold. First, a preliminary result is needed.

\textbf{Theorem 3.} Assume that for  $f_k$, Properties 1 and 3 hold  and let $(q(t): t \in T_k)$ be a distribution such that
$$
sgn \left( f_k (q(t): t \in T_k) \right) =  0.
$$
Then, for all distributions $(p(t): t \in T_k)$ ,
$$
sgn \left( f_k (p(t)+q(t): t \in T_k)  \right) =  sgn \left( f_k(p(t): t \in T_k)  \right).
$$

\textbf{Proof.}
The distributions $(q(t): t \in T_k)$ and $(p(t): t \in T_k)$  may be seen as distributions in two layers of a $k+1$-dimensional table, of which $T_k$ is the marginal table.

If $ f_k (p(t): t \in T_k)  $ is zero, directional collapsibility implies the result immediately.

If $ f_k (p(t): t \in T_k)$  is positive, then $(p(t)+q(t): t \in T_k)$ will be written as the sum of two distributions, so that $f_k$ is positive on both, and, then, it is also positive on $(p(t)+q(t): t \in T_k)$. 

Because of  Property 1, the entry in $p(1, 1 , \ldots 1)$  may be decreased by a positive amount, such that the entry remains positive and  $  f_k (p(t): t \in T_k)$    also remains positive. If the entry  $q(1, 1 , \ldots 1)$  is increased by the same amount, then by (i), $ f_k (q(t): t \in T_k) $  becomes positive, and by directional collapsibility,  $f_k(p(t)+q(t): t \in T_k)$ has to be positive.

If $ f_k (p(t): t \in T_k)$ is negative, the argument is modified so that  $p(1, 1 , \ldots 1)$  is increased by a small amount. 

\hfill{ }\qed

\textbf{Theorem 4.} If for a  parameter of association $f_k$, Properties 1 and 2 hold, then Property 3 holds  for it if and only if, for any distribution,
$$
sgn \left( f_k (p(t): t \in T_k)\right) = sgn \left( DI_k(f_k (p(t): t \in T_k)\right) . 
$$

\textbf{Proof.}  
The ''if'' part follows from the directional collapsibility of $DI_k$.
The idea of the proof of the "only if'' part is to write $ (p(t): t \in T_k)$ as the sum of several $k$-dimensional distributions, such that $f_k$ has, on at least one of them, the same sign as $DI_k$ does, and on the others it either has the same sign or is zero. Then, repeated application of directional collapsibility and of Theorem 3 yields the result.

Let $s$ be the smallest entry in $(p(t), t \in T_k)$. Subtract from every entry in  $(p(t), t \in T_k)$ the value $s$, to obtain $(p(t)-s, t \in T_k)$, which has non-negative entries  and will be denoted as $q_0(t)$. Then
\begin{equation}\label{sameDI}
DI_k(p(t): t \in T_k) = DI_k(q_0(t): t \in T_k).  
\end{equation}

Let $(u_0(t), t \in T_k)$ be zero in every cell of $T_k$. If the following condition does not hold for $l=0$,
\begin{equation} \label{emptyhalf}
q_l(t)=0 \makebox { for all } t \in T_{ke} \makebox{ or } q_l(t)=0 \makebox { for all } t \in T_{ke},
\end{equation}
then consider the smallest positive entry (or one of the smallest positive ones), of $(q_0(t), t \in T_k)$. Suppose it is in cell $t_1$. The cell $t_1$ is  either in $T_{ke}$ or in $T_{ko}$, and there is a cell, say $t'_1$ in the other subset, so that $q_0(t'_1) \ge q_0(t_1)$. Then define
$$u_1(t) = q_0(t_1), \makebox{ if } t=t_1 \makebox{ or } t=t'_1, 
$$
$$
u_1(t)=0 \makebox{ otherwise,}
$$
and set
$$
q_1(t) = q_0(t)-u_1(t), \,\, t \in T_k.
$$

Continue the procedure of the previous paragraph for $q_1$ instead of $q_0$ to obtain $u_2$ and $q_2$, and repeat until the condition in (\ref{emptyhalf}) wil become true for some $l$.

Now the entries of  $q_l$ are all zero, either in $T_{ko}$ or in  $T_{ke}$ or in both, and because 
$$
\sum_{j=0}^l \sum_{t \in T_{ke}} u_j(t) = \sum_{j=0}^l \sum_{t \in T_{ko}} u_j(t) ,
$$   
this will happen when $DI_k(p(t): t\in T_k)$ is positive or negative or zero, respectively.

In case $DI_k(p(t): t \in T_k) > 0$, define for every $t' \in T_{ke}$, and in case $DI_k(p(t): t \in T_k) < 0$, define for every $t' \in T_{ko}$
$$
v_{t'}(t') = q_l(t'), 
$$
$$
v_{t'}(t'') = 0 \makebox{ if } t'' \neq t'.
$$

It follows from the construction that if $DI_k(p(t): t \in T_k) > 0 $, then 
$$
p(t) - s = q_0(t) = \sum_{0=1}^l u_j(t) + \sum_{t' \in T_{ke}} v_{t'}(t),
$$
if  $DI_k(p(t): t \in T_k) < 0$, then
$$
p(t) - s = q_0(t) = \sum_{0=1}^l u_j(t) + \sum_{t' \in T_{ko}} v_{t'}(t),
$$
and if $DI_k (p(t): t \in T_k)= 0$, then
$$
p(t) - s = q_0(t) = \sum_{0=1}^l u_j(t) ,
$$
for every $t \in T_k$.

In the first two cases, add to all the enties in te distributions $(u_j, j = 0, \ldots , l)$, and to all the  entries in the distributions $(v_{t'}, t' \in T_{ke})$ (or $(v_{t'}, t' \in T_{ko})$), the value of
$$
\frac{s}{(l+1) 2^k + 2^{k-1}2^k},
$$ 
and in the third case, to all the  enties in the distributions $(u_j, j = 0, \ldots , l)$, the value of
$$
\frac{s}{(l+1)2^k},
$$
and denote the distributions obtained by  $u'_j, j = 0, \ldots , l$ and  $v'_{t'}, t' \in T_{ke}$ (or $v'_{t'}, t' \in T_{ko}$), so that the following holds:

If $DI_k(p(t): t \in T_k) > 0 $, then 
\begin{equation} \label{pos}
p(t)  = \sum_{0=1}^l u'_j(t) + \sum_{t' \in T_{ke}} v'_{t'}(t),
\end{equation}
if  $DI_k(p(t): t \in T_k) < 0$, then
\begin{equation} \label{neg}
p(t)  = \sum_{0=1}^l u'_j(t) + \sum_{t' \in T_{ko}} v'_{t'}(t),
\end{equation}
and if $DI_k(p(t): t \in T_k)= 0$, then
\begin{equation} \label{zeroo}
p(t)  = \sum_{0=1}^l u'_j(t),
\end{equation}
for every $t \in T_k$.

The structures of the distributions in (\ref{pos}),  (\ref{neg})  and  (\ref{zeroo}) are as follows:

The distribution $u'_0$ has the same entry in every cell.

Each of the distributions $u'_j, j = 1, \ldots, l$ has the same entry in every cell, except for one cell in $T_{ke}$ and one cell in $T_{ko}$, which have the same value in them (different from the other cells) and these will be called the specific entries.

Each of the distributions $v'_{t'}, t' \in T_{ke}$ ($v'_{t'}, t' \in T_{ko}$) has the same entry in every cell, which will be called the common value, except for a cell in $T_{ke}$  ( $T_{ko}$), which has a larger value.

To complete the proof, it will be shown that 

\begin{equation} \label{fu}
f_k(u'_j) = 0, \,\, j = 0, 1, \ldots l 
\end{equation}
\begin{equation} \label{fpos}
f_k(v'_{t'}) > 0, \,\, t' \in T_{ke}
\end{equation}
 \begin{equation} \label{fneg}
f_k(v'_{t'}) < 0, \,\, t' \in T_{ko}
\end{equation}
which, together with directional collapsibility, Theroem 3, and (\ref{pos}),  (\ref{neg})  and  (\ref{zeroo}), imply the desired result.

To see (\ref{fu}) for $j=0$, note that because all entries are the same, swapping the categories of one variable does not change the distribution but changes the sign of $f_k$ to its opposite by (ii), thus $f_k(u'_0)=0$.

To see (\ref{fu}) for $j=1, \ldots, l$, consider a series of swaps of indices of variables, which exchanges the two cells with the specific values. If such a series of swaps exists, it leaves the distribution unchanged, as all other entries are the same. Such a series of swaps is obtained, if the indices of all variables are swapped, in an arbitrary order, which are $2$ in any one of the cells with specific values but not in the other one.  One of these cells is in $T_{ke}$, thus has an even number of $2'$s, the other cell is in  $T_{ko}$, thus has an odd number of $2'$s. Therefore, the total number of indices equal to $2$ in the two cells is odd. To obtain the number of indices that are equal to $2$ in exactly one of the cells, from the odd total, the number of $2'$s in identical positions in the indices has to be subtracted. This latter number is even, so the total number of swaps is odd. By repeated application (ii), the sign of $f_k$ changes to its opposite during the series of swaps, but because the distribution remains the same, it cannot change. Thus, $f_k$ is zero for $u'_j, j=1, \ldots \l$.     

To see (\ref{fpos}) and  (\ref{fneg}), note first that if $t'=(1, 1, \ldots, 1)$, then $f_k$ is positive, because, as was seen in the proof of  (\ref{fu}) for $j=0$, for a distribution with all entries equal, $f_k$ is zero, and if the entry in the cell $(1, 1, \ldots, 1)$ is increased to the value  $v'_{(1, 1, \ldots, 1)}(1, 1, \ldots, 1)$, then by (i), $f_k$ will become positive. 

For any  $t' \in T_{k}$, other than $(1, 1, \ldots, 1)$, write the value of $v'_{t'}(t')$ in celll (1, 1, \ldots, 1), while keeping the common value in the other cells. This will give a positive $f_k$. This entry can be moved into the cell $t'$ by a series of swaps, while the common value remains in all other cells and also appears in $(1, 1, \ldots, 1)$. This requires an even number of swaps, if $t' \in T_{ke}$, keeping the positive value, thus (\ref{fpos})  is implied,  and an odd number of swaps, if   
  $t' \in T_{ko}$, yielding a negative value, thus (\ref{fneg})  is implied.

\hfill{ }\qed

\section{Discussion}

This section addresses briefly the meaning and use of the results of the paper.

Odds ratios and log-linear parameters have the very attractive property of being variation independent from lower dimensional marginals, and, thus, make it possible to identify  association with the information in the joint distribution which is there {\em in addition} to the information in the lower dimensional marginal distribution, see \cite{rudas98}. In particular, Property 4 implies that variants of the Iterative Proportional Fitting / Scaling algorithm may be used to obtain maximum likelihood estimates in various exponential family models that are specified by prescribing the values of odds ratios \citep{bfh75, rudas91, kr14}. However, as implied by Theorem 2, this property  makes it impossible to find parameters of association, which are free from the possibility of Simpson's paradox, if Properties 1 and 2 are assumed. 

The lack of directional collapsibility is considered problematic by most analysts, as testified not only by a large body of literature about ''avoiding'' it, but also by the wide-spread use of the Mantel-Haensel odds ratio in meta-analysis, which always estimates the common odds ratio to be  in between the lowest and highest conditional odds ratios, even if the marginal odds ratio is outside of this range.

On the other hand, as implied by Theorem 4, if only the direction of association is of interest, and one wishes to use parameters of association which are directionally collapsible, then, if Properties 1 and 2 are assumed, there is only one possible choice for this direction, and it is given by $DI_k$. 

The simple linear contrast of the cell probabilities,  $DI_k$, is not necessarily seen as a meaningful  parameter of association, and those who are not willing to accept the direction of association as given by it, have to accept that  Simpson's paradox cannot be avoided. Another argument for using $DI_k$ in certain situations, given in \cite{rudas10}, is that if the data are observational, then allocation in treatment categories is potentially informative, thus association (effect) should not be measured by a parameter which is variationally independent of the treatment marginal(s). In such cases, avoiding Simpson's paradox is an additional bonus, which comes with using $DI_k$ 

Whether or not one is ready to adopt $DI_k$ to determine the direction of association, it is worth noting that its sampling behaviour is straightforward, in particular it does not depend on the individual cell entries, and not even on $k$. If the population probability (fraction) of cells in $T_{ke}$ is $p$, then the probability that $sgn(DI_k)=1$, which will lead to the correct or incorrect decision as to the direction of association depending on whether $p>0.5$ or $p \le 0.5$, may be obtained as follows. In the case of multinomial sampling with $N$ observations,
$$
P(DI_k > 0) =  \sum_{x=[N/2]+1}^{N}  \binom{N}{x} p^x (1-p)^{N-x},
$$ 
Which, for large sample sizes, may be approximated as
$$
 \Phi \left(  \sqrt{N} \frac{p-0.5}{\sqrt{p(1-p)}}   \right),
$$
where $\Phi$ is  the cumulative distribution function of the standard normal distribution. For example, with a sample size of $1000$, and true value of $DI_k= 0.05$, that is $p=0.525,$ the probability of correctly deciding that the association is positive is about $0.94$, which seems quite certain, even though the assumed true value is not very far from zero. An important property of the probability of correct decision with $DI_k$, is that it (in addition to the sample size), only depends on the true value of $DI_k$. In contrast, the probability of correct decision with $LOR_k$, depends, in addition to the sample size, also on the individual cell probabilities.

\section{Acknowledgment}

 This research was supported in
part by Grant K-106154 from the Hungarian National Scientiﬁc Research Fund (OTKA). The author is indebted to Anna Klimova and Ren\'ata N\'emeth for several helpful comments.

\end{document}